\documentclass[preprint,12pt]{elsarticle}

\usepackage{amssymb}

\usepackage{amsmath}
  \usepackage{paralist}
  \usepackage{graphics} 
  \usepackage{epsfig} 
 \usepackage[colorlinks=true]{hyperref}
\usepackage{ifpdf}

\allowdisplaybreaks[4]

\usepackage{graphicx,amssymb,mathrsfs,latexsym,amsthm}
\usepackage{verbatim}

\newtheorem{theorem}{Theorem}
\newtheorem{lemma}[theorem]{Lemma}
\newtheorem{corollary}[theorem]{Corollary}

\newtheorem{definition}[theorem]{Definition}
\newtheorem{example}[theorem]{Example}

\begin{document}

\begin{frontmatter}



\title{Projection Pressure and Bowen's Equation for a Class of Self-similar Fractals with Overlap Structure}


\author[Nanjing]{Chenwei Wang}
\ead{chenweiwang01@163.com}
\author[Nanjing,NJU]{Ercai Chen\fnref{NSFC}}
\ead{ecchen@njnu.edu.cn}
\address[Nanjing]{School of Mathematics, Nanjing Normal
University\\Nanjing 210097, China}
\address[NJU]{Center of Nonlinear Science, Nanjing University\\Nanjing 210093,
China}

\fntext[NSFC]{The author is supported by the National Natural
Science Foundation of China (10971100) and National Basic Research
Program of China (973 Program) (2007CB814800).}

\begin{abstract}
Let $\{S_i\}_{i=1}^{l}$ be an iterated function system(IFS) on
$\mathbb{R}^d$ with attractor K. Let $\pi$ be the canonical
projection. In this paper we define a new concept called
``projection pressure" $P_\pi(\varphi)$ for $\varphi\in
C(\mathbb{R}^d)$ under certain affine IFS, and show the variational
principle about the projection pressure. Furthermore we check that
the unique zero root of ``projection pressure'' still satisfies
Bowen's equation when each $S_i$ is the similar map with the same
compression ratio. Using the root of Bowen's equation, we can get
the Hausdorff dimension of the attractor $K$.
\end{abstract}

\end{frontmatter}


\section{Introduction}

  \indent \indent Let $\{S_i:X \rightarrow X\}_{i=1}^l$ be a family of contractive maps on a nonempty closed
  set $X\subset \mathbb{R}^d$. Following Barnsley \cite{[Bar]}, we say that $\{S_i\}_{i=1}^l$ is
  an \emph{iterated function system} (IFS) on $X$. Hutchinson \cite{[Hutchinson]} showed that there is a unique
  nonempty compact set $K \subset X$, called the \emph{attractor} of $\{S_i\}_{i=1}^l$, such that
  $K=\cup_{i=1}^lS_i(K)$.

  There are many references to compute the Hausdorff dimension of $K$ or the Hausdorff dimension of multifractal
  spectrum, such as \cite{[McMullen]}, \cite{[Moorthy]} and \cite{[King]}. Thermodynamic formalism played a significant role when we try to compute the Hausdorff dimension of $K$, especially the Bowen's equation. Usually we call
  $P_J(t\psi)=0$ the Bowen's equation, where $P_J$ is the topological pressure of the map $f:J \rightarrow J$, and $\psi$ is the geometric potential $\psi(z)=\log|f'(z)|$. The root of Bowen's equation always approaches the Hausdorff dimension of some sets. In \cite{[Bowen]}, Bowen first discovered this equation while studying the Hausdorff dimension of quasi-circles. Later Ruelle \cite{[Ruelle]}, Gatzouras and Peres \cite{[Gatzouras]} showed that Bowen's equation gives the Hausdorff dimension of $J$ whenever $f$ is a $C^1$ conformal map on a Riemannian manifold and $J$ is a repeller.
   According to the method for calculating Hausdorff dimension of cookie-cutters presented by Bedford \cite{[Bedford1]}, Keller discussed the relation between classical pressure and dimension for IFS \cite{[Keller]}. He concluded that
  if $\{S_i\}_{i=1}^l$ is a conformal IFS satisfying the disjointness condition with local energy function $\psi$, then the
  pressure function has a unique zero root $t_0=\dim_H K$. In 2000, using the definition of Carath$\acute{e}$odory dimension characteristics, Barreira and Schmeling \cite{[BS]} introduced the notion of the $u$-dimension $\dim_uZ$ for positive functions $u$, showing that $\dim_uZ$ is the unique number $t$ such that $P_Z(-tu)=0$.

   On the progress of calculating $\dim_HK$, references \cite{[McMullen]}, \cite{[Moorthy]} and \cite{[King]} depend on the
  open set condition and  separable condition. In fact, there are a lot of examples don't satisfy this disjointness
  condition. Hui Rao and Zhiying Wen once discussed a kind of self-similar fractal with overlap structure called $\lambda$-Cantor set
  \cite{[Rao]}.

  In order to study the Hausdorff dimension of an invariant measure $\mu$ for
  conformal and affine IFS with overlaps, Dejun Feng and Huyi Hu introduce a notion \emph{projection
  entropy} (see\cite{[Feng]}), which plays the similar role as the classical entropy of IFS satisfying the open set condition,
  and it becomes the classical entropy if the projection is finite to one.

  Bedford pointed out that the Bowen's equation works if
  three elements are present:(i) conformal contractions, (ii) open set
  conditions, and (iii) subshift of finite type (Markov) structure. Chen
  \cite{[chen2]} proved that subshift of finite type (Markov)
  structure can be replaced by any subshift structure. The motivation of this paper is to find a new pressure function
satisfying the Bowen's equation without the open set conditions. We
will show the variational principle for the new pressure under
certain affine IFS, and check that the zero of projection pressure
is equal to $\dim_H K$ when each $S_i$ has the same compression
ratio. Thus, using our results we can easily get the dimension of
${\lambda}$ cantor set.

 \section{The projection pressure for certain affine IFS: Definition and variational principle}

 \indent \indent Let $\{S_i\}_{i=1}^l$ be an IFS on a closed set $X\subset \mathbb{R}^d$. Denote by $K$ its attractor.
 Let $\sum=\{1, \cdots, l\}^\mathbb{N}$ associated with the left shift $\sigma$. Let $M_\sigma(\Sigma)$ denote the space
 of $\sigma$-invariant measure on $\Sigma$, endowed with the weak-star topology, $C(X)$ the space of real-valued continuous functions of $X$ and $\pi:\Sigma\rightarrow K$ be the
 canonical projection defined by
 \begin{equation}\label{1}
    \{\pi(x)\}=\bigcap_{n=1}^{\infty}S_{x_1} \circ S_{x_2} \circ \cdots \circ S_{x_n}(K),\quad \textrm { where }
    \quad x=\{x_i\}_{i=1}^{\infty}.
 \end{equation}

  A measure $\mu$ on $K$ is called invariant (resp.ergodic) for the IFS if there is an invariant (resp.ergodic)
  measure $\nu$ on $\Sigma$ such that $\mu=\nu\circ\pi^{-1}$.

 Let ($\Omega, \;\mathcal{F},\; \nu$) be a probability space. For a sub-$\sigma$-algebra $\mathcal{A}$ of $\mathcal{F}$
 and $f\in L^1(\Omega, \;\mathcal{F},\; \nu)$, we denote by $E_{\nu}(f|\mathcal{A})$ the conditional expectation of $f$
 given $\mathcal{A}$. For countable $\mathcal{F}$-measurable partition $\xi$ of $\Omega$. We denote by
 $I_{\nu}(\xi|\mathcal {A})$ the conditional information of $\xi$ given $\mathcal{A}$, which is given by the formula
 \begin{equation}\label{2}
    I_{\nu}(\xi|\mathcal{A})=-\sum_{A\in {\xi}} \mathcal {X}_A  \lg{E_\nu {(\mathcal {X}_A|\mathcal{A})}},
 \end{equation}
where  $\mathcal {X}_A$ denotes the characteristic function on
$\mathcal{A}$.

The conditional entropy of $\xi$ given $\mathcal{A}$, written
$H_{\nu}(\xi|\mathcal{A})$ is defined by the formula
$H_{\nu}(\xi|\mathcal{A})=\int I_{\nu}(\xi|\mathcal{A})d\nu$.

 The above information and entropy are unconditional when $\mathcal{A}=\mathcal{N}$, the trivial $\sigma$-algebra
 consisting of sets of measure zero and one, and in this case we write
 \begin{equation}\label{3}
    I_{\nu}(\xi|\mathcal{A})=I_\nu(\xi)   \quad \textrm{and} \quad H_{\nu}(\xi|\mathcal{N})=H_\nu(\xi).
 \end{equation}

 Now we consider the space$(\Sigma, \;\mathscr{B}(\Sigma), \;m)$, where $\mathscr{B}(\Sigma)$ is the Borel
 $\sigma$-algebra on $\Sigma$ and $m\in M_\sigma(\Sigma)$. Let $\mathcal {P}$ denote the Borel partition
 \begin{equation*}
    \mathcal {P} =\{[j]:1 \leq j \leq l \}
 \end{equation*}
  of $\Sigma$, where $[j]=\{(x_i)_{i=1}^\infty \in \Sigma, \; x_1=j\}$. Let $\mathcal {I}$ denote the $\sigma$-algebra
 \begin{equation*}
    \mathcal{I}=\{B\in\mathcal{B}(\Sigma):\sigma^{-1}B=B\}.
 \end{equation*}

 For convenience, we use $\gamma$ to denote the Borel $\sigma$-algebra
 $\mathcal {B}(\mathbb {R}^d)$ of $\mathbb {R}^d$. For $f\in C(X)$, denote $\|f\|=\sup_{x\in X} f(x)$ and $S_nf(x)=\sum_{i=0}^{n-1}f(\sigma^{n} x)$  $x\in X$.
Let $\Sigma_n=\{[b]: [b]=(x_1, x_2,\cdots,x_n), \; x_i \in \Sigma,\;
i=1,\cdots,n\} $.

\begin{definition}\label{projection entropy}
For any $m\in
M_\sigma(\Sigma)$, we call
\begin{equation*}\label{projection entropy}
    h_\pi(\sigma,m)=H_m(\mathcal{P}|\sigma^{-1}\pi^{-1}\gamma)-H_m(\mathcal{P}|\pi^{-1}\gamma)
\end{equation*}
the \emph{projection entropy} of $m$ under $\pi$ \emph{w.r.t.}
$\{S_i\}_{i=1}^l$ , and we call
\begin{equation*}
    h_\pi(\sigma,m,x)=E_m(f|\mathcal {I})(x)
\end{equation*}
the \emph{local projection entropy of m at $x$ under $\pi$ w.r.t
}$\{S_i\}_{i=1}^l$, where $f$ denotes the function
$I_m(\mathcal{P}|\sigma^{-1}\pi^{-1}\gamma)-I_m(\mathcal{P}|\sigma^{-1}\gamma)$.
\end{definition}

It is clear that $h_{\pi}(\sigma,m)=\int h_{\pi}(\sigma,m,x)dm(x)$.

\begin{definition}
Let $k\in \mathbb{N}$ and $\nu \in M_{\sigma^k}(\Sigma)$.
Define
\begin{equation*}
    h_{\pi}(\sigma^k,\nu)=H_\nu(\mathcal{P}_0^{k-1}|\sigma^{-k}\pi^{-1} \gamma)-H_\nu(\mathcal{P}_0^{k-1}|\pi^{-1}\gamma).
\end{equation*}
\end{definition}

The term $h_{\pi}(\sigma^k,\nu)$ can be viewed as the projection
measure-theoretic entropy of $\nu$ w.r.t. the IFS$\{ S_{i_1}\circ
\cdots \circ S_{i_k}: 1\leq i_j \leq l \, \rm{for} \, 1\leq j \leq
k\}.$ The following lemma exploits the connection between
$h_\pi(\sigma^k, \nu)$ and $h_\pi(\sigma, \nu)$, where
$m=\frac{1}{k}\sum_{i=0}^{k-1}\nu \circ \sigma^{-i}$.

\begin{lemma}\label{lem1}
Let $k\in \mathbb{N}$ and $\nu \in M_{\sigma^k}{\Sigma}.$ Set
$m=\frac{1}{k}\sum_{i=0}^{k-1}\nu \circ \sigma^{-i}$. Then $m$ is
$\sigma$-invariant, and $h_\pi(\sigma,
\nu)=\frac{1}{k}h_\pi(\sigma^k, \nu)=\frac{1}{k}h_\pi(\sigma^k, m)$.
\end{lemma}
{\textbf{Proof}}: See Proposition 4.3 in\cite{[Feng]}.\\

The following two Lemmas will be used in our results.
\begin{lemma}\label{lem2}
\textrm{Let} $a_1, a_2,\cdots, a_k$ \textrm{be given real numbers.
If }$p_i\geq0$ \textrm{and} $\sum_{i=0}^kp_i=1$, then
$\sum_{i=0}^kp_i(a_i-\log{p_i})\leq \log{(\sum_{i=0}^ke^{a_i})}$
\textrm{and equality holds iff}
$p_i=\frac{e^{a_i}}{\sum_{j=1}^ke^{a_j}}$.
\end{lemma}
{\textbf{Proof}}: See Lemma 9.9 in \cite{[Walters]}.\\

\begin{lemma}\label{lem3}
Assume that $\Omega$ is a subset of $\{1,\cdots,l\}$ such that
$S_i(K)\bigcap S_j(K)=\emptyset$ for all $i,j\in \Omega $ with
$i\neq j$. Suppose that $\nu$ is an invariant measure on $\Sigma$
supposed on $\Omega^{\mathbb{N}}$, i.e., $\nu([j])=0$ for all
$j\in\{1,\cdots,l\}\backslash \Omega$. Then
$h_\pi(\sigma,\nu)=h(\sigma,\nu)$.
\end{lemma}
{\textbf{Proof}}: See Lemma 4.19 in \cite{[Feng]}.\\

Let $S_i(x)=Ax+c_i   \quad  (i=1, \cdots l)$ be an IFS on $\mathbb{R}^d$, where $A$ is a $d\times d$ non-singular matrix with $\|A\|<1$ and $c_i\in \mathbb{R}^d$. Let $K$ denote its attractor and $\pi: \Sigma \rightarrow K$ be the canonical projection. Let $\mathcal{Q}$ denote the partition $\{[0,1)^d+\alpha: \alpha \in \mathbb{Z}^d\}$ of $\mathbb{Z}^d$. We set $\mathcal{Q}_n=\{A^nQ:Q\in \mathcal{Q}\}$ for $n=0,1, \cdots$.\\

\begin{lemma}\label{lem4}
Let $m\in M_\sigma(\Sigma)$, then
$h_\pi(\sigma,m)=\lim\limits_{n\rightarrow\infty}
\frac{H_m(\pi^{-1}\mathcal{Q}_n)}{n}$.
\end{lemma}
{\textbf{Proof}}: See Proposition 4.18 (i) in \cite{[Feng]}.\\

\begin{theorem}
 If an IFS $\{S_i\}_{i=1}^l$ has the form as
above and $f \in C(K)$.
 Then

\[ \lim_{n \rightarrow \infty}  \frac{1}{n}(\log \sum_{\substack { Q\in
\mathcal{Q}_n \\  Q \cap K\neq \emptyset} } \sup_{\alpha \in
\pi^{-1}(Q \cap K)}e^{S_nf\pi(\alpha)})=\sup\{h_\pi(\sigma,m)+\int
f\circ \pi dm: \quad  m \in M_\sigma (\Sigma)  \}. \]
\end{theorem}

{\textbf{Proof}}: We assume $K\subset[0,1)^d$, without loss of
generality.
We divided the proof into two steps.\\                         
$step1$. \[ \liminf_{n\rightarrow\infty} \frac{1}{n} (\log
\sum_{\substack { Q \in\mathcal{Q}_n  \\   Q\cap K\neq \emptyset } }
\sup_{\alpha \in \pi^{-1}(Q\cap K)}e^{S_nf\pi(\alpha)})\geq \sup
\{h_\pi(\sigma,m)+\int f\circ \pi dm: \quad  m \in M_\sigma(\Sigma)
\}  .                                           \] For arbitrary $n
\in \mathbb{N}, Q \in \mathcal{Q}_n, m \in M_{\sigma}(\Sigma)$, let
$g_n(Q)=\sup_{\alpha \in \pi^{-1}(Q \cap K)}S_n f \pi (\alpha),
P(Q)=m(\pi^{-1}Q)$. By lemma \ref{lem2}, we have
\begin{eqnarray*}
\log{\sum _{\substack{Q \in \mathcal{Q}_n \\ Q\cap K \neq
\emptyset}} \sup_{\alpha \in \pi^{-1}(Q \cap K)}e^{S_nf\pi
(\alpha)}}
&\geq& \sum_{\substack{Q \in \mathcal{Q}_n \\ Q\cap K \neq \emptyset}}P(Q)(g_n(Q)-\log P(Q)) \\
&=&  H_m(\pi^{-1} \mathcal{Q}_n)+\sum_{\substack{Q \in \mathcal{Q}_n \\ Q \cap K \neq \emptyset}}P(Q)g_n(Q)\\
&\geq&   H_m(\pi^{-1} \mathcal{Q}_n)+\int S_nf\pi (\alpha)dm \\
&=&  H_m(\pi^{-1} \mathcal{Q}_n)+n \int f \pi dm.
\end{eqnarray*}
Using Lemma \ref{lem4} yields
\begin{eqnarray*}
\liminf_{n\rightarrow \infty} \frac{1}{n}(\log \sum_{\substack{Q \in
\mathcal{Q}_n \\ Q\cap K \neq \emptyset}} \sup_{\alpha \in
\pi^{-1}(Q\cap K)}e^{S_nf\pi (\alpha)})
&\geq& \lim_{n\rightarrow \infty}\frac{H_m(\pi^{-1}\mathcal{Q}_n)}{n}+\int f\pi dm\\
&=&h_\pi(\sigma,m)+\int f\pi dm.
\end{eqnarray*}
By the arbitrariness of $m$, we have $step1$.\\            
$step2.$ \[   \sup \{h_\pi(\sigma,m)+\int f\circ \pi dm: \quad  m
\in M_\sigma(\Sigma)\} \geq  \limsup_{n\rightarrow\infty}
\frac{1}{n} (\log \sum_{\substack { Q \in\mathcal{Q}_n  \\   Q\cap
K\neq \emptyset } } \sup_{\alpha \in \pi^{-1}(Q\cap K)
}e^{S_nf\pi(\alpha)}).                                           \]
 By the continuity of $f\pi$, for arbitrary $\epsilon > 0$, there exists $N \in \mathbb{N}$ such that for arbitrary
 $a_N\in \Sigma_N$, and any $x,y\in a_N$ we have\\
 \begin{equation*}
    |f\pi(x)-f\pi(y)|<\epsilon.
 \end{equation*}
 Now, for any $n\in \mathbb{N}$ and $a_{n+N}\in \Sigma_{n+N}$
 \begin{equation*}
    |S_{n+N}f\pi(x)-S_{n+N}f\pi(y)|\leq n\epsilon+2N\|f\pi\| \qquad \forall x,y \in a_{n+N}.
 \end{equation*}
 For $x\in K$ and $n\in \mathbb{N}$, we denote $$g_n(x)=\sup_{\alpha \in \pi^{-1}x}e^{S_nf\pi(\alpha)}.$$
 For $X\subset \mathbb{R}^d$ and $X\cap K\neq \emptyset$, we denote $$g_n(X)=\sup_{\substack{x\in X\cap K}}g_n(x).$$
 For any $Q\in \mathcal{Q}_n$, we denote $$2Q=\bigcup_{\substack{ P \in \mathcal{Q}_n\\ \overline{P}\cap \overline{Q}
 \neq \emptyset}}P$$ and $$3Q=\bigcup_{\substack{ P \in \mathcal{Q}_n\\ \overline{P}\cap \overline{2Q} \neq \emptyset}}P.$$
 Claim :There exists a subset $\Gamma$ of $\{ {Q \in \mathcal{Q}_{n+N}:Q\cap K \neq \emptyset} \}$ such that

  (i) $\sum\limits_{Q \in \Gamma} g_{n+N}(Q) \geq \frac{\sum\limits_{Q\in \mathcal {Q}_{n+N}}g_{n+N}(Q)}{7^d}$;

  (ii) $2Q\cap 2\widetilde{Q}=\emptyset ,  \quad Q ,\ \widetilde{Q}\in \Gamma ,\quad Q\neq \widetilde{Q}$.\\

 Since $K$ is compact, we can find $Q_1\in \mathcal{Q}_{n+N}$ and $x_1\in \mathcal{Q}_1\cap K $
 such that $$g_{n+N}(x_1)=g_{n+N}(Q\cap K)=g_{n+N}(K).$$
 If $K\setminus 3Q=\emptyset$, then $\Gamma=\{Q_1\}.$
 Otherwise, let $K_2=K\setminus 3Q_1^{\circ} $ we can find $Q_2\in \mathcal{Q}_{n+N}$ and $x_2\in Q_2\cap K$ such that $$g_{n+N}(x_2)=
 g_{n+N}(Q_2\cap K)=g_{n+N}(K).$$ If $K_2\setminus3Q_2=\emptyset$, then $\Gamma=\{Q_1, Q_2\}$.
 Repeat above steps and we can finish it in finite steps.
 Clearly, there exists $M=M(n,N)$ and $\Gamma = \{Q_1, Q_2,\cdots,Q_M \} $
 such that $\rm{(i),(ii)}$ are satisfied. Let $X=\{x_1, x_2, \cdots, x_M\},$
 according to our claim, we have
 \begin{equation}\label{4}
    \sum_{x \in X}g_{n+N}(x)\geq \frac{\sum\limits_{\substack{Q\in \mathcal{Q}_{n+N} \\ Q \cap K \neq \emptyset}}
    g_{n+N}(Q)}{7^d}.
 \end{equation}
 For each $x\in X$, since $x \in K$, we can pick a word $[u]_x\in \Sigma_{n+N}$
 such that $x\in S_{[u]_x}K$ and
 \begin{equation}\label{3.5}
    \sup_{\alpha \in [u]_x}e^{(S_{n+N}f\pi)(\alpha)}\geq \sup_{\alpha \in \pi^{-1}x}e^{(S_{n+N}f\pi)(\alpha)}.
\end{equation}
Consider the collection $W=\{[u]_x, \,x\in X\}$. The separation
condition for elements in $X$ guarantees that $S_{[u]_x}K \cap
S_{[u]_y}K= \emptyset$ for all $x,y \in X$ with $x\neq y.$ Define a
Bernoulli measure $\nu$ on $W^{\mathbb{N}}$ by
\begin{equation*}
    \nu([u]_x)=\frac{\sup\limits_{\alpha \in [u]_x}e^{(S_{n+N}f\pi)(\alpha)}}{\sum\limits_{x\in X}\sup\limits_{\beta \in
    [u]_x}e^{(S_{n+N}f\pi)(\beta)}},
\end{equation*}
\begin{equation*}
    \nu([w_1,\cdots,w_k])=\prod_{i=1}^k\nu([w_i]), \quad  w_i\in W, \;k\in\mathbb{N}.
\end{equation*}
Then $\nu$ can be viewed as a $\sigma^{n+N}$-invariant measure on
$\Sigma$ (by viewing $W^\mathbb{N}$ as a subset of $\Sigma$). By
Lemma \ref{lem3}, we have $h_\pi(\sigma^{n+N}, \nu)=h(\sigma^{n+N},
\nu).$
Define $\mu= \frac{1}{n+N} \sum\limits_{i=0}^{n+N-1}\nu\circ\sigma ^{-i}\in M_\sigma(\Sigma)$\\
and $\xi=\{[u]_x,\,x\in X\} \cup \{\Sigma\setminus \bigcup_{x\in
X}[u]_x\}$.

 According to lemma \ref{lem1}, we have
\begin{eqnarray*}
  h_\pi(\sigma, \mu)+\int f\pi d\mu
   &=& \frac{h_\pi ((\sigma^{n+N},\nu))}{n+N}+\frac{\int S_{n+N} f \pi d\nu}{n+N} \\
  &=&\frac{1}{n+N}(h(\sigma^{n+N},\nu)+\int S_{n+N}f\pi d\nu)\\
  &=&\frac{1}{n+N}(H_v(\xi)+\int S_{n+N}f\pi d\nu)\\
  &\geq& \frac{1}{n+N}\bigg(\sum_{x\in X}\Big(-\nu([u]_x)\log \nu([u]_x)+\nu([u]_x)\inf_{\alpha \in [u]_x}S_{n+N}f\pi(\alpha)\Big)\bigg)\\
  &\geq& \frac{1}{n+N}\bigg(\sum_{x\in X}\Big(-\nu([u]_x)\log \nu([u]_x)+ \nu([u]_x)(\sup_{\alpha \in [u]_x}S_{n+N}f\pi(\alpha)\\
  &&-n\epsilon-2N\|f\pi\|)\Big)\bigg)\\
  &=&\frac{1}{n+N}\bigg(\sum_{x\in X}\Big(-\nu([u]_x)\log \nu([u]_x)+\nu([u]_x) \sup_{\alpha \in [u]_x}S_{n+N}f\pi(\alpha)\Big)\bigg)\\
  &&-\frac{n\epsilon+2N\|f\pi\|}{n+N}\\
  &=&\frac{1}{n+N}\log {\sum_{x \in X}\sup_{\alpha \in [u]_x}e^{(S_{n+N}f\pi)(\alpha)}}- \frac{n\epsilon+2N\|f\pi\|}{n+N}\\
  &\geq& \frac{1}{n+N}\log {\sum_{x \in X}\sup_{\alpha \in \pi^{-1}x}e^{(S_{n+N}f\pi)(\alpha)}}- \frac{n\epsilon+2N\|f\pi\|}{n+N}\\
  &\geq& \frac{1}{n+N}\log {\bigg(\frac {\sum\limits_{\substack{Q \in \mathcal{Q} \\Q\cap K \neq \emptyset}}
g_{n+N}(Q)}{7^d}\bigg)}-\frac{n\epsilon+2N\|f\pi\|}{n+N}.\\
  \end{eqnarray*}

Let $k=n+N$ and let $n\rightarrow\infty$, then $k\rightarrow\infty$.
We have
\begin{equation*}
    \sup \{h_\pi(\sigma,m)+\int f\circ \pi dm , m \in M_\sigma(\Sigma)\} \geq  \limsup_{n\rightarrow\infty}\frac{1}{n}
    (\log \sum_{\substack { Q \in\mathcal{Q}_n  \\   Q\cap K\neq \emptyset } } \sup_{\alpha \in \pi^{-1}(Q\cap K)}
    e^{S_nf\pi(\alpha)})-\epsilon.
\end{equation*}
Since $\epsilon$ is arbitrary, we finish the proof of $step2$.

\begin{definition}\label{projection pressure}
If $f\in C(\mathbb{R}^d, \mathbb{R})$,
we call
\begin{equation*}\label{projection pressure}
    P_\pi(f)=\lim_{n \rightarrow \infty}  \frac{1}{n}(\log \sum_{\substack { Q\in \mathcal{Q}_n \\
    Q \cap K\neq \emptyset} } \sup_{\alpha \in \pi^{-1}(Q \cap K)}e^{S_nf\pi(\alpha)})
\end{equation*}
the projection pressure of $f$ under $\pi$ w.r.t.$\{S_i\}_{i=1}^l$,
where $\{S_i\}_{i=1}^l$ as in Theorem 2.1.
\end{definition}

It is clearly that if $f=0$ we have the same result of Proposition
4.18 (ii) in \cite{[Feng]}.
\begin{corollary} $\lim\limits_{n\rightarrow \infty} \frac{\log \#\{Q\in \mathcal{Q}:A^nQ \cap K \neq \emptyset  \}}{n}
=\sup\{h_{\pi}(\sigma,m):\, m \in M_{\sigma}(\Sigma)\}$.
\end{corollary}

\section{Bowen's equation for certain self-similar set}
\begin{definition}\label{conformal}
The IFS $\{S_i\}_{i=1}^l$ is conformal if for every $i\in
\{1,2,\cdots,l\}$, (1) $S_i: U \rightarrow S_i(U)$ is $C^1$, (2)
$\|S'_i(x)\| \neq 0$ for all $x\in U$ , (3)
$|S_i'(x)y|=\|S_i'(x)\||y|$ for all $x\in U$, $y \in \mathbb{R}^d$.
\end{definition}

In the following, we assume that $S_i(x)=Ax+c_i   \quad  (i=1,
\cdots l)$ be an IFS and $A$ be a $d\times d$ compressed orthogonal
matrix, which means $A$ satisfies $AA^T=cE$, $c_i \in \mathbb{R}^d$
and $0<c<1$. Clearly, such a IFS is conformal. Let $\|S\|$ denote
the spectral norm of $S$.

\begin{lemma}\label{3.1} 
Let $K$ be the attractor of a conformal IFS $S_i(x)=Ax+c_i$
$(i=1,\cdots,l)$. Then we have
\begin{eqnarray}
  \dim_HK  &=& \frac{\sup \{h_\pi(\sigma,m): m\in M_\sigma(\Sigma)\}}{-\log\|A\|} .
\end{eqnarray}
\end{lemma}
{\textbf{Proof}}: See Theorem 2.13 in \cite{[Feng]}.\\

\begin{theorem} Let $S_i(x)=Ax+c_i   \quad  (i=1, \cdots l)$ be
an IFS and $A$ be a $d\times d$ compressed orthogonal matrix. Let
$\pi: \Sigma \rightarrow K$ be the canonical projection. Then
$\dim_H K$ is the unique root of $P_\pi(\log\|A\| \cdot t)=0$.
\end{theorem}

{\textbf{Proof}}: According to Theorem 2.1 and Lemma \ref{3.1}, we have \\
\begin{equation*}
    P_\pi(\log\|A\|\cdot t)=\sup\{h_\pi(\sigma,m): \quad m\in M_\sigma(\Sigma)\}+\log \|A\| \cdot t.
\end{equation*}
Hence $\dim_H K$ is the unique root of $P_\pi(\log\|A\| \cdot t)=0$.\\

Thus we verify the projection pressure function satisfies the
Bowen's equation without any disjoint property.

\begin{example}$(\lambda -Cantor \, \,set)$
 Let $\lambda \in [0,1]$ be a real number and let
$S_1(x)=x/3,S_2(x)=x/3+\lambda/3,S_3(x)=x/3+2/3$ be three
similarities on $\mathbb{R}$. Then the self-similar set generated by
these three similarities, denoted by $F_\lambda$, will be called a
$\lambda$-Cantor set.

If $\lambda$=0, then $S_1=S_2$ and $F_\lambda$ is exactly the
classical middle third Cantor set. In this case
$\dim_H(K)=\log2/\log3$.

If $\lambda$=1, then $F_\lambda=[0,1]$ and $\dim_H(K)=1$.

If $0<\lambda<1$, the structure of $F_\lambda$ is quite complicated.

However, according to our result, the Hausdorff dimension of
$F_\lambda$ should be the unique root $t_0$ of $P_{\pi}(-\log 3\cdot
t)=0$.
\end{example}

\begin{figure}[h!]
\centering
\ifpdf
  \setlength{\unitlength}{1bp}%
  \begin{picture}(406.74, 181.52)(0,0)
  \put(0,0){\includegraphics{11.ps}}
  \put(182.83,18.73){\fontsize{14.23}{17.07}\selectfont x}
  \put(12.86,163.24){\fontsize{14.23}{17.07}\selectfont y}
  \put(12.02,24.75){\fontsize{14.23}{17.07}\selectfont O}
  \put(393.98,19.01){\fontsize{14.23}{17.07}\selectfont x}
  \put(215.27,160.38){\fontsize{14.23}{17.07}\selectfont y}
  \put(215.80,23.27){\fontsize{14.23}{17.07}\selectfont O}
  \end{picture}%
\else
  \setlength{\unitlength}{1bp}%
  \begin{picture}(406.74, 181.52)(0,0)
  \put(0,0){\includegraphics{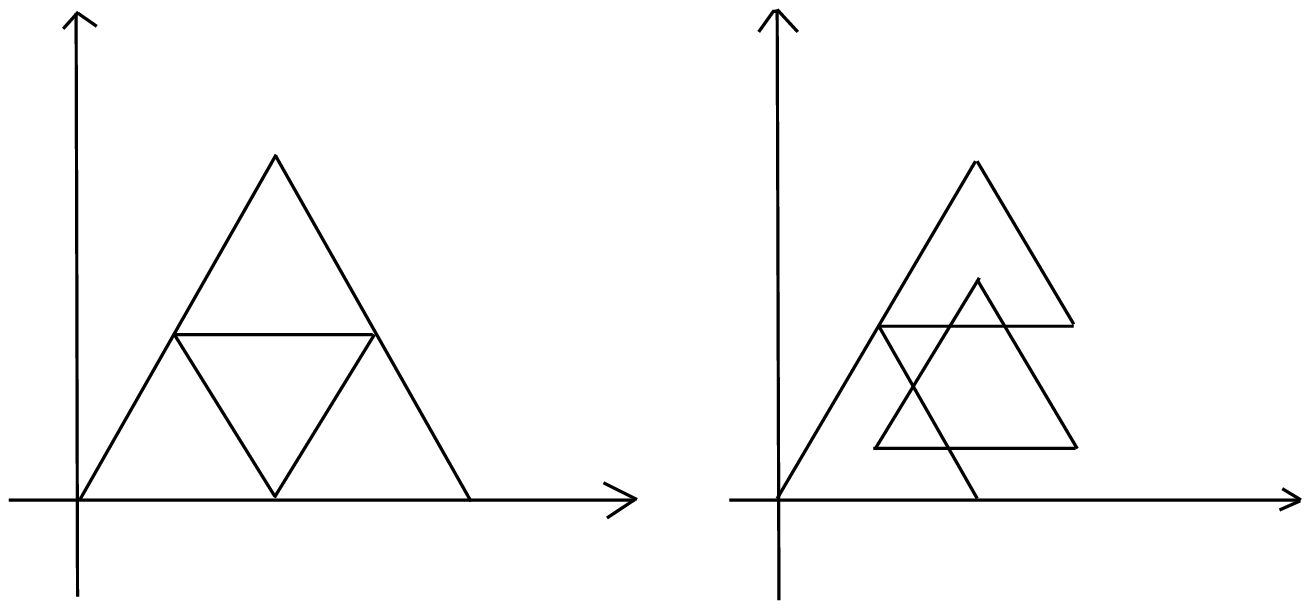}}
  \put(182.83,18.73){\fontsize{14.23}{17.07}\selectfont x}
  \put(12.86,163.24){\fontsize{14.23}{17.07}\selectfont y}
  \put(12.02,24.75){\fontsize{14.23}{17.07}\selectfont O}
  \put(393.98,19.01){\fontsize{14.23}{17.07}\selectfont x}
  \put(215.27,160.38){\fontsize{14.23}{17.07}\selectfont y}
  \put(215.80,23.27){\fontsize{14.23}{17.07}\selectfont O}
  \end{picture}%
\fi \caption{\protect\label{fig:1}}
\end{figure}

\begin{example}$(\textrm{Overlapping Sierpi$\acute{n}$ski triangle})$
Consider a certain IFS as :

$S_1\left(
                                                \begin{array}{c}
                                                  x \\
                                                  y \\
                                                \end{array}
                                              \right)
=\left(
     \begin{array}{cc}
       1/2 & 0 \\
       0 & 1/2 \\
     \end{array}
   \right)
\left(
  \begin{array}{c}
    x \\
    y \\
  \end{array}
\right) $, $S_2\left(
          \begin{array}{c}
            x \\
            y \\
          \end{array}
        \right)=\left(
                  \begin{array}{cc}
                    1/2 & 0 \\
                    0 & 1/2 \\
                  \end{array}
                \right)\left(
                         \begin{array}{c}
                           x \\
                           y \\
                         \end{array}
                       \right)+\left(
                                 \begin{array}{c}
                                   a_1 \\
                                   a_2 \\
                                 \end{array}
                               \right)
$, where $0\leq a_1 \leq 1/2$ and $0 \leq a_2 \leq \sqrt{3}/4$,
$S_3\left(
          \begin{array}{c}
            x \\
            y \\
          \end{array}
        \right)=\left(
                  \begin{array}{cc}
                    1/2 & 0 \\
                    0 & 1/2 \\
                  \end{array}
                \right)\left(
                         \begin{array}{c}
                           x \\
                           y \\
                         \end{array}
                       \right)+\left(
                                 \begin{array}{c}
                                   1/4 \\
                                   \sqrt{3}/4 \\
                                 \end{array}
                               \right)
$ on $\{(x,y):y\leq \sqrt{3}x,y \geq 0 \,\rm{and}\, y\geq
-\sqrt{3}x+\sqrt{3} \}$(See Figure\ref{fig:1}). Let $K_{a_1, a_2}$
be its attractor.

If $a_1=a_2=0$ or $a_1=1/4$ and $a_2=\sqrt{3}/4$, then $\dim_H
K_{a_1, a_2}=1$.

If $a_1=1/2$ and $a_2=0$, then $K_{a_1, a_2}$ is classical
Sierpi$\acute{n}$ski triangle and its Hausdorff dimension is $\log 3
/ \log 2$.

If $0< a_1 < 1/2$ and $0 < a_2 < \sqrt{3}/4$,  the Hausdorff
dimension of $K_{a_1, a_2}$ should be the unique root $t_0$ of
$P_{\pi}(-\log 2\cdot t)=0$.
\end{example}

\begin{thebibliography}{99}
\bibitem{[Bar]} M. Barnsley, Fractals everywhere, Academic Press, Inc., Boston, MA, (1988).
\bibitem{[BS]} L. Barreira and J. Schmeling, Sets of ``non-typical" points have full topological entropy and full Hausdorff dimension, Israel J. Math., 116 (2000), 29-70.
\bibitem{[Bedford1]} T. Bedford, Applications of dynamical systems to fractals - a study of cookie-cutter Cantor sets, In Fractal Geometry and Analysis, Kluwer, (1991), 1-44.
\bibitem{[Bedford2]} T. Bedford, Hausdorff dimension and box dimension in self-similar sets, Topology and Measure V, University of Greifswald, (1987).
\bibitem{[Bowen]} R. Bowen, Hausdorff dimension of quasi-circles, Publ. Math. Inst. Hautes $\acute{\textrm{E}}$tudes Sci., 50 (1979), 259-273.
\bibitem{[Chen1]} E. C. Chen and J. C. Xiong, Dimension and measure theoretic entropy of a subshift in symbolic space, Chinese Sci. Bull., 42 (1997), 1193-1196.
\bibitem{[chen2]} E. C. Chen and Y. S. Sun, The Bowen's formula for Hausdorff dimensions of invariant sets, Preprint.
\bibitem{[Feng]} D. J. Feng and H. Y. Hu, Dimension theory of iterated function systems, Comm. Pure Appl. Math., 62 (2009), 1435-1500.
\bibitem{[Gatzouras]} D. Gatzouras and Y. Peres, Invariant measures of full dimension for some expanding maps, Ergodic Theory Dynam. Systems, 17(1) (1997), 147-167.
\bibitem{[Hutchinson]} J. E. Hutchinson, Fractals and self-similarity, Indiana Univ. Math. J., 30 (1981), 713-747.
\bibitem{[Keller]} G. Keller, Equilibrium states in ergodic theory, Cambridge University Press, Cambridge, (1998).
\bibitem{[King]} J. King, The singularity spectrum for general Sierpi$\acute{\textrm{n}}$ski carpets, Adv. Math., 116 (1995), 1-8.
\bibitem{[Mane]} R. Ma$\tilde{\textrm{n}}\acute{\textrm{e}}$, Ergodic theory and differentiable dynamics, Springer-Verlag, Berlin, (1987).
\bibitem{[McMullen]} C. McMullen, The Hausdorff dimension of general Sierpi$\acute{\textrm{n}}$ski carpets, Nagoya Math. J., 96 (1984), 1-9.
\bibitem{[Moorthy]} C. G. Moorthy, R. Vijaya and P. Venkatachalapathy, Hausdorff dimension of Cantor-like sets, Kyungpook Math. J., 32 (1992), 197-202.
\bibitem{[Pesin]} Y. Pesin, Dimension theory in dynamical systems, contemporary views and applications, University of Chicago Press, (1998).
\bibitem{[Rao]} H. Rao and Z. Y. Wen, A class of self-similar fractals with overlap structure, Adv. in Appl. Math., 20 (1998), 50-72.
\bibitem{[Rohlin]} V. A. Rohlin, On the fundamental ideas of measure theory, Mat. Sbornik N.S., 25(67) (1949), 107-150.
\bibitem{[Rudin]} W. Rudin, Real and complex analysis, Third edition, McGraw-Hill Book Co., New York, (1987).
\bibitem{[Ruelle]} D. Ruelle, Repellers for real analytic maps, Ergodic Theory Dynam. Systems, 2(1) (1982), 99-107.
\bibitem{[Stoll]} R. R. Stoll, Linear algebra and matrix theory, McGraw-Hill Company, Inc., New York-Toronto-London, (1952).
\bibitem{[Walters]} P. Walters, An introduction to ergodic theory, Graduate Texts in Mathematics, 79, Springer, Berlin, (2000).
\end{thebibliography}
\end{document}